\documentclass[10pt]{article}
\usepackage{amsfonts}
\usepackage{amsmath}
\usepackage{mathrsfs}
\usepackage{mathrsfs,amscd,amssymb,amsthm,amsmath,bm,graphicx,psfrag,subfigure}

\makeatletter

\renewcommand{\@seccntformat}[1]{{\csname the#1\endcsname}{\normalsize .}\hspace{.5em}}
\makeatother

\def \[{\begin{equation}}
\def \]{\end{equation}}

\newtheorem{thm}{Theorem}[section]
\newtheorem{prop}{Proposition}
\newtheorem{defi}{Definition}

\newtheorem{lem}[thm]{Lemma}
\newtheorem{cor}[thm]{Corollary}

\newtheorem{prob}[thm]{Problem}

\newenvironment{wst}
{\setlength{\leftmargini}{1.5\parindent}
 \begin{itemize}
 \setlength{\itemsep}{-1.1mm}}
{\end{itemize}}

\begin{document}
\setlength{\baselineskip}{15pt}
\begin{center}{\Large \bf On the spectral moment of trees with given degree sequences\footnote{Financially supported by the National Natural
Science Foundation of China (Grant Nos. 11071096,\, 11271149) and the Special Fund for Basic Scientific Research of Central Colleges (CCNU11A02015).}}

\vspace{4mm}

{\large Shuchao Li\footnote{E-mail: lscmath@mail.ccnu.edu.cn (S.C.
Li), 425333559@qq.com (Y.B. Song)},\ \ Yibing Song}\vspace{2mm}

Faculty of Mathematics and Statistics,  Central China Normal
University, Wuhan 430079, P.R. China
\end{center}
\noindent {\bf Abstract}: Let $A(G)$ be the adjacency matrix of graph $G$ with eigenvalues $\lambda_1(G), \lambda_2(G),\dots, \lambda_n(G)$ in non-increasing order. The number $S_k(G):=\sum_{i=1}^{n}\lambda_i^{k}(G)\, (k=0, 1,\ldots, n-1)$ is called the $k$th spectral moment of $G$. Let $S(G) = (S_0(G), S_1(G),\ldots, S_{n-1}(G))$ be the sequence of spectral moments of $G.$ For two graphs $G_1, G_2$, we have $G_1\prec_{s}G_2$ if for some $k \in \{1,2,3,\dots,n-1\}$, we have $S_i(G_1) = S_i(G_2)\, ,\, i = 0, 1,\ldots, k-1$ and $S_k(G_1)<S_k(G_2).$ In this paper, the last $n$-vertex tree with a given degree sequence in an $S$-order is determined. Consequently, we also obtain the last trees in an $S$-order in the sets of all trees of order
$n$ with the largest degree, the leaves number, the independence number and the matching number, respectively.

\vspace{2mm} \noindent{\it Keywords}: Spectral moment; Tree; Degree sequence

\vspace{2mm}


 {\setcounter{section}{0}
\section{\normalsize Introduction}\setcounter{equation}{0}
All graphs considered here are finite, simple and connected. Undefined terminology and notation may be referred to \cite{D-I}.
Let $G=(V_G,E_G)$ be a simple undirected graph with $n$ vertices. $G-v$, $G-uv$ denote the graph obtained from $G$ by deleting vertex $v \in V_G$, or edge
$uv \in E_G$, respectively (this notation is naturally extended if more than one vertex or edge is deleted). Similarly,
$G+uv$ is obtained from $G$ by adding an edge $uv \not\in E_G$. For $v\in V_G,$ let
$N_G(v)$ (or $N(v)$ for short) denote the set of all the adjacent vertices of $v$ in $G$ and $d_G(v)=|N_G(v)|$. A \textit{leaf} of a graph is a vertex of degree one.

Let $A(G)$ be the adjacency matrix of a graph $G$ with $\lambda_1(G),\lambda_2(G),\dots,\lambda_n(G)$ being its eigenvalues in non-increasing order. The number $\sum_{i=1}^n\lambda_i^k(G)\, (k=0,1,\dots,n-1)$ is called the $k$th spectral moment of $G$, denoted by $S_k(G)$. Let $S(G)=(S_0(G), S_1(G), \dots, S_{n-1}(G))$ be the sequence of spectral moments of $G$. For two graphs $G_1, G_2$, we shall write $G_1=_sG_2$ if $S_i(G_1)=S_i(G_2)$ for $i=0,1,\dots,n-1$. Similarly, we have $G_1\preceq_sG_2\, (G_1$ comes before $G_2$  in an $S$-order) if for some $k\, (1\leqslant k\leqslant {n-1})$, we have  $S_i(G_1)=S_i(G_2)\, (i=0,1,\dots,k-1)$  and $S_k(G_1)<S_k(G_2)$. We shall also write $G_1\preceq_sG_2$ if $G_1\prec_sG_2$ or $G_1=_sG_2$.  $S$-order has been used in producing graph catalogs (see \cite{C-G}), and for a more general setting of spectral moments one may be referred to \cite{C-R-S1}.

Investigation on $S$-order of graphs attracts more and more researchers' attention. Cvetkovi\'c and Rowlinson \cite{C-R-S3} studied
the $S$-order of trees and unicyclic graphs and characterized the first and the last graphs, in an $S$-order,
of all trees and all unicyclic graph with given girth, respectively. Chen, Liu and Liu \cite{C-L-L} studied the lexicographic ordering
by spectral moments ($S$-order) of unicyclic graph with a given girth. Wu and Fan \cite{D-F} determined the first
and the last graphs, in an $S$-order, of all unicyclic graphs and bicyclic graphs, respectively. Pan et al. \cite{X-F-P}
gave the first $\sum_{k=1}^{\lfloor\frac{n-1}{3}\rfloor}(\lfloor\frac{n-k-1}{2}\rfloor-k+1)$ graphs apart from an $n$-vertex path, in an $S$-order, of all trees with $n$ vertices. Wu and Liu \cite{D-M} determined the last $\lfloor\frac{d}{2}\rfloor+1$, in an $S$-order,
among all $n$-vertex trees of diameter $d\, (4 \le d \le n-3)$. Pan et al. \cite{B-Z1} identified the
last and the second last graphs, in an $S$-order, of quasi-trees. Hu and Li \cite{H-Li} studied the spectral moments of graphs with given number of clique number and chromatic number, respectively. Li and Zhang \cite{Li-Z} studied the spectral moments of graphs with given number of cut edges.

A nonincreasing sequence of nonnegative integers $\pi = (d_0, d_1, . . . , d_{n-1})$ is called \textit{graphic} if there exists a graph
having $\pi$ as its vertex degree sequence. Motivated by the recent results in terms of vertex degrees, we generally propose
the following question.
\begin{prob}
For a given graphic degree sequence $\pi$, let
$$
  \text{$\mathscr{G}_{\pi} = \{G|G$ is connected with $\pi$ as its degree sequence\}.}
$$
Characterize the last (first) graph in an $S$-order among all graphs $G$ in $\mathscr{G}_{\pi}$.
\end{prob}

In this paper, we only consider a special case for the above problem, i.e., for a given degree sequence of some tree. The main result of this paper is as follows:
\begin{thm}
For a given degree sequence $\pi$ of some $n$-vertex tree, let
$$
  \text{$\mathscr{T}_{\pi} = \{T|\, T$ is an $n$-vertex tree with $\pi$ as its degree sequence\}.}
$$
Then $T^*$ (see in Section 2) is a unique last tree in an $S$-order, among $\mathscr{T}_\pi$.
\end{thm}

The rest of the paper is organized as follow. In Section 2, some notations and preliminary results are presented. In Section 3, we present the proof of Theorem 1.2 and some corollaries.

\section{\normalsize Preliminary }\setcounter{equation}{0}
Throughout we denote by $P_n, K_{1,n-1}$ the path, star on $n$ vertices, respectively. Let $B_5$ be a tree obtained from $P_3$ by attaching two pendant vertices to one of its end vertices; it is easy to see that the degree sequence of $B_5$ is $(3,2,1,1,1)$.

\begin{lem} [\cite{B-Z1}]
The $k$th spectral moment of $G$ is equal to the number of closed walks of length $k$.
\end{lem}

Let $F$ be a graph. An $F$-subgraph of $G$ is a subgraph of $G$ which is isomorphic to the graph $F.$ Let $\phi_G(F)$ (or $\phi(F)$ for short) be the number of all $F$-subgraphs of $G$.

\begin{lem}
For every graph $G$, we have
\begin{wst}
\item[{\rm (i)}]$ S_4(G)=2\phi(P_2)+4\phi(P_3)+8\phi(C_4)$ {\rm (see \cite{D-M})};
\item[{\rm (ii)}]$ S_5(G)=30\phi(C_3)+10\phi(U_4)+10\phi(C_5)$ {\rm(see \cite{C-R-S1})}.
\end{wst}
\end{lem}

Since all the graphs we considered are trees, note that $S_i(T_1) = S_i(T_2)$ for $i=0,1,2,3,4,5,7$ where $T_1,T_2 \in \mathscr{T}_\pi.$ Moreover, by Lemma 1.1 we can get
\begin{eqnarray}
  S_6(T_1)-S_6(T_2) &=& 6(\phi_{T_1}(P_4)-\phi_{T_2}(P_4)), \\
  S_8(T_1)-S_8(T_2) &=& 32(\phi_{T_1}(P_4)-\phi_{T_2}(P_4))+8(\phi_{T_1}(P_5)-\phi_{T_2}(P_5))+16(\phi_{T_1}(B_5)-\phi_{T_2}(B_5)).
\end{eqnarray}

For a given non-increasing degree sequence $\pi = (d_0,d_1,\dots,d_{n-1})$ of a tree with $n \geqslant 3$, we use breadth-first search method to define a special tree $T^*$ with degree sequence $\pi$ as follows (see also \cite{X-D-Z}). Assume that $d_m > 1$ and $d_{m+1}= \cdots = d_{n-1}=1$ for $0 \leqslant m < n-2$. Put $s_0 = 0.$ Select a vertex $v_{01}$ as a root and begin with $v_{01}$ in layer 0. Put $s_1 = d_0$ and select $s_1$ vertices $\{v_{11},\dots,v_{1,s_1}\}$ in layer 1 such that they are adjacent to $v_{01}.$ Thus $d(v_{01})=s_1=d_0.$ We continue to construct all other layers by recursion. In general, put $s_t=d_{s_0+s_1+ \cdots +s_{t-2}+1} + \cdots + d_{s_0+s_1+ \cdots +s_{t-2}+s_{t-1}}-s_{t-1}$ for $t \geqslant 2$ and assume that all vertices in layer $t$ have been constructed and are denoted by $\{v_{t,1},\dots,v_{t, s_t} \}$ with $d(v_{t-1,1})=d_{s_0+ \cdots +s_{t-2}+1},\dots,d(v_{t-1,s_{t-1}})=d_{s_0+\cdots+s_{t-1}}.$ Now using the induction hypothesis, we construct all vertices in layer $t+1$. Put $s_{t+1}=d_{s_0+\dots+s_{t-1}+1}+\cdots+d_{s_0+\cdots+s_t}-s_t.$ Select $s_{t+1}$ vertices $\{v_{t+1,1},\ldots,v_{t+1,s_{t+1}} \}$ in layer $t+1$ such that $v_{t+1,i}$ is adjacent to $v_{tr}$ for $r=1$ and $1 \leqslant i \leqslant d_{s_0+\dots+s_{t-1}+1}-1$ and for $2 \leqslant r \leqslant s_t$ and $d_{s_0+\dots+s_{t-1}+1}+d_{s_0+\dots+s_{t-1}+2}+\cdots+d_{s_0+\dots+s_{t-1}+r-1}-r+2 \leqslant i \leqslant d_{s_0+\dots+s_{t-1}+1}+d_{s_0+\dots+s_{t-1}+2}+\cdots+d_{s_0+\dots+s_{t-1}+r}-r.$ Thus $d(v_{tr})=d_{s_0+\dots+s_{t-1}+r}$ for $1 \leqslant r \leqslant s_t.$ Assume that $m=s_0+\cdots+s_{p-1}+q.$ Put $s_{p+1}=d_{s_0+\dots+s_{p-1}+1}+\cdots+d_{s_0+\dots+s_{p-1}+q}-q$ and select $s_{p+1}$ vertices $\{v_{p+1,1},\dots,v_{p+1,s_{p+1}} \}$ in layer $p+1$ such that $v_{p+1,i}$ is adjacent to $v_{pr}$ for $1 \leqslant r \leqslant q$ and $d_{s_0+\dots+s_{p-1}+1}+d_{s_0+\cdots+s_{p-1}+2}+\cdots+d_{s_0+\dots+s_{p-1}+r-1}-r+2 \leqslant i \leqslant d_{s_0+\dots+s_{p-1}+1}+d_{s_0+\dots+s_{p-1}+2}+\cdots+d_{s_0+\dots+s_{p-1}+r}-r.$ Thus $d(v_{p,i})=d_{s_0+\dots+s_{p-1}+i}$ for $1 \leqslant i \leqslant q.$ In this way, we obtain a tree $T^*.$ It is easy to see that $T^*$ is of order $n$ with degree sequence $\pi$.
\begin{figure}[h!]
\begin{center}
  \psfrag{a}{$v_{01}$}\psfrag{e}{$v_{14}$}
  \psfrag{b}{$v_{11}$}\psfrag{f}{$v_{21}$}
  \psfrag{c}{$v_{12}$}\psfrag{d}{$v_{13}$}
  \psfrag{n}{$v_{29}$}\psfrag{l}{$v_{27}$}
  \psfrag{k}{$v_{26}$}\psfrag{i}{$v_{24}$}
  \psfrag{j}{$v_{25}$}\psfrag{g}{$v_{22}$}
  \psfrag{h}{$v_{23}$}\psfrag{m}{$v_{28}$}
  \psfrag{o}{$v_{31}$}\psfrag{p}{$v_{32}$}
  \psfrag{q}{$v_{33}$}\psfrag{2}{$T$}\psfrag{t}{$v_{2,10}$}
  \psfrag{1}{$T^*$}
  \includegraphics[width=160mm]{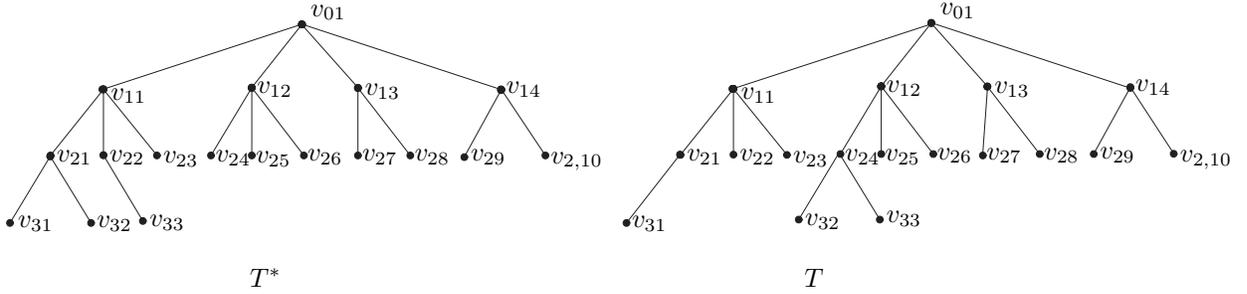}\\
  \caption{Trees $T^*$ and $T$ with the same degree sequences $\pi=(4,4,4,3,3,3,2,1,1,1,1,1,1,1,1,1,1,1).$}
\end{center}
\end{figure}

For example, for a given degree sequence $\pi=(4,4,4,3,3,3,2,1,1,1,1,1,1,1,1,1,1,1),$ $T^*$ is the tree of order 18; see Fig. 1. There is a vertex $v_{01}$ in layer 0; four vertices $v_{11},v_{12},v_{13},v_{14}$ in layer 1; ten vertices $v_{21},v_{22},\dots,v_{29}, v_{2,10}$ in layer 2; three vertices $v_{31},v_{32},v_{33}$ in layer 3. Moreover, $s_0=0,s_1=d_0=4,s_2=d_1+d_2+d_3+d_4-s_1=4+4+3+3-4=10,$ $s_3=d_5+\cdots+d_{13}-s_2=3,m=s_1+q=4+2=6.$

For a graph with a root $v_0$, we call the distance the $\textit{height}$ $h(v)=dist(v,v_0)$ of a vertex $v.$

\begin{defi}
Let $T=(V,E)$ be a tree with root $v_0.$ A well-ordering $\prec$ of the vertices is called breadth-first search ordering with non-increasing degrees ($BFS$-ordering for short) if the following holds for all vertices $u,v \in V$:
\begin{wst}
\item[{\rm (1)}] $u \prec v $ implies $h(u) \leqslant h(v);$

\item[{\rm (2)}] $u \prec v $ implies $d(u) \geqslant d(v);$

\item[{\rm (3)}] if there are two edges $uu_1 \in E(T)$ and $vv_1 \in E(T)$ such that $u \prec v,h(u)=h(u_1)+1$ and $h(v)=h(v_1)+1,$ then $u_1 \prec v_1.$
\end{wst}
\end{defi}

We call trees that have a BFS-ordering of its vertices a BFS-tree. All trees have an ordering which satisfy the conditions (1) and (3) by using breadth-first search, but not all trees have a BFS-ordering. For example, the tree $T$ (see Fig. 1) of order 18 has not a $BFS$-ordering with degree sequence $\pi=(4,4,4,3,3,3,2,1,1,1,1,1,1,1,1,1,1,1).$ In fact, Zhang obtained the following result in 2008.
\begin{prop}[\cite{X-D-Z}]
For a given degree sequence $\pi$ of some tree, there exists a unique tree $T^*$ with degree sequence $\pi$ having a $BFS$-ordering. Moreover, any two trees with the same degree sequences and having $BFS$-ordering are isomorphic.
\end{prop}
We recall the notion of majorization. Let $\pi=(d_0,\dots,d_{n-1})$ and $\pi_{'}$=$(d_0',\dots,d_{n-1}')$ be two non-increasing sequences. If $\sum_{i=0}^kd_i \leqslant \sum_{i=0}^kd_i'$ and $\sum_{i=0}^{n-1}d_i = \sum_{i=0}^{n-1}d_i'$, then the sequence $\pi'$ is said to \textit{major} the sequence $\pi$ and denoted by $\pi \lhd \pi'$. It is known that the following result holds.
\begin{prop}[\cite{P-E-T}]
 Let $\pi=(d_0,\dots,d_{n-1})$ and $\pi'=(d_0',\dots,d_{n-1}')$ be two non-increasing graphic degree sequences. If $\pi \lhd \pi',$ then there exists a series graphic degree sequences $\pi_1,\dots,\pi_k$ such that $\pi \lhd \pi_1 \lhd \cdots \lhd \pi_k \lhd \pi',$ and only two components of $\pi_i$ and $\pi_{i+1}$ are different from $1$.
\end{prop}

\section{\normalsize Main results}\setcounter{equation}{0}

In order to show Theorem 1.2, we need the following lemmas. The first lemma follows immediately by its definitions.
\begin{lem}
For any $n$-vertex tree $T$, one has
\begin{eqnarray}
\phi_T(P_4)&=&\sum_{uv \in E_T}(d_T(u)-1)(d_T(v)-1).\label{eq:3.1}\\
\phi_T(T_5) &=& \sum_{uv \in E_T}\left[{d_T(u)-1 \choose 2}(d_T(v)-1)+{d_T(v)-1 \choose 2}(d_T(u)-1)\right].\label{eq:3.2}\\
\phi_T(P_5) &=& \sum_{u \in V_T}\sum_{x, y\in N_T(u)}(d_T(x)-1)(d_T(y)-1).\label{eq:3.3}
\end{eqnarray}
\end{lem}
\begin{lem}
Let $T=(V_T, E_T)$ be a tree with $v_1u_1, v_2u_2 \in E_T,$ and $v_1v_2,\, u_1u_2 \notin E_T.$ Let $T'=T-\{u_1v_1, u_2v_2\}+\{u_1u_2, v_1v_2\}.$ If $d_T(v_1)> d_T(u_2)$ and $d_T(v_2)> d_T(u_1),$ then $T \prec_s T'.$
\end{lem}
\begin{proof}
By Lemma 1.1,\, $S_i(T)=S_i(T')$ for $i=0, 1, 2, 3, 5$. Note that $T$ and $T'$ have the same degree sequence, hence $\phi_T(P_3)=\sum_{x\in V_T}{d_T(x) \choose 2}=\phi_{T'}(P_3)$. By Lemma 2.2(i), we obtain that $S_4(T)=S_4(T')$.
Combining with (2.1) and (3.1) we have
\begin{eqnarray*}
S_6(T')-S_6(T)&=&  6(\phi_{T'}(P_4)-\phi_{T}(P_4))\\
               &=& 6[(d_{T}(v_1)-1)(d_{T}(v_2)-1)+(d_{T}(u_1)-1)(d_{T}(u_2)-1)\\
                &&-(d_{T'}(v_1)-1)(d_{T'}(u_1)-1)+(d_{T'}(v_2)-1)(d_{T'}(u_2)-1)]\\
               &=&6(d_{T}(v_1)-d_{T}(u_2))(d_{T}(v_2)-d_{T}(u_1))\\
               &>& 0,
\end{eqnarray*}
which implies that $T \prec_s T'$.
\end{proof}

In view of the proof of Lemma 3.1, we may let $\sigma=\max\{\phi_T(P_4): T\in \mathscr{T}_{\pi}\}.$ Set \text{$\mathscr{T}'_{\pi}:=\{T: T\in \mathscr{T}_{\pi}$ with} $\phi_T(P_4)=\sigma\}$. Hence, the last tree $T^{*}$ must be contained in $\mathscr{T}_\pi'.$

\begin{lem}
Let $T$ be an arbitrary tree rooted at one maximum degree vertex. Assume that $u, v$ in $V_T$ satisfy $d_T(u)=d_T(v)$ and $\sum_{x \in N_T(u)}d_T(x) \geqslant \sum_{x \in N_T(v)}d_T(x)$.  Set
\begin{eqnarray*}
  d_T(x_0)&:=&\text{$\min\{d_T(x)|x \in N_T(u)$ with $h(x)=h(u)+1\}$},\\
  d_T(x_1)&:=&\text{$\max\{d_T(x)|x \in N_T(v)$ with $h(x)=h(v)+1\}$}.
\end{eqnarray*}
Let $T'=T-\{ux_0, vx_1\}+\{ux_1, vx_0\}$. If $d_T(x_0)<d_T(x_1)$, then $T\prec_s T'.$
\end{lem}
\begin{proof}
By Lemma 1.1,\, $S_i(T)=S_i(T')$ for $i=0, 1, 2, 3, 5, 7$. Note that $T$ and $T'$ have the same degree sequence, hence $\phi_T(P_3)=\sum_{x\in V_T}{d_T(x) \choose 2}=\phi_{T'}(P_3)$. By Lemma 2.2(i), we obtain that $S_4(T)=S_4(T')$. In view of (2.1) and (\ref{eq:3.1}), we have
\begin{eqnarray*}
S_6(T')-S_6(T)&=&  6(\phi_{T'}(P_4)-\phi_{T}(P_4))\\
               &=& 6[(d_{T'}(u)-1)(d_{T'}(x_1)-1)+(d_{T'}(v)-1)(d_{T'}(x_0)-1)\\
                &&-(d_T(u)-1)(d_T(x_0)-1)+(d_T(v)-1)(d_T(x_1)-1)]\\
               &=&0.
\end{eqnarray*}
The last equality follows from $d_T(u)=d_{T'}(u)=d_T(v)=d_{T'}(v),\, d_T(x_0)=d_{T'}(x_0)$  and $d_T(x_1)=d_{T'}(x_1).$  Hence, $S_6(T')=S_6(T)$ and
\[\label{eq:3.4}
 \phi_{T'}(P_4)=\phi_{T}(P_4).
\]

Notice that $T$ and $T'$ have the same degree sequences, in view of (\ref{eq:3.2}) we have
\[\label{eq:3.5}
\phi_T(T_5)=\phi_{T'}(T_5).
\]

In view of (2.2) and (\ref{eq:3.3})-(\ref{eq:3.5}) we have
\begin{eqnarray*}
S_8(T')-S_8(T)&=&  8(\phi_{T'}(P_5)-\phi_{T}(P_5))\\
              & =& 8\left[\sum_{u \in V_{T'}}\sum_{x,y\in N_{T'}(u)}(d_{T'}(x)-1)(d_{T'}(y)-1) -
               \sum_{u \in V_T}\sum_{x,y\in N_T(u)}(d_T(x)-1)(d_T(y)-1)\right]\\
               &=&8\left[\sum_{x \in N_{T'}(u)\backslash x_1}(d_{T'}(x_1)-1)(d_{T'}(x)-1)+\sum_{x \in N_{T'}(v)\backslash x_0}(d_{T'}(x_0)-1)(d_{T'}(x)-1)\right.\\
               &&\left.-\sum_{x \in N_{T}(u)\backslash x_0}(d_{T}(x_0)-1)(d_{T}(x)-1)-\sum_{x \in N_{T}(v)\backslash x_1}(d_{T}(x_1)-1)(d_{T}(x)-1)\right]\\
               &=&8(d_T(x_1)-d_T(x_0))\left[\sum_{x \in N_{T}(u)\backslash x_0}(d_{T}(x)-1)-\sum_{x \in N_{T}(v)\backslash x_1}(d_{T}(x)-1)\right]\\
               &=& 8(d_T(x_1)-d_T(x_0))\left[\left(\sum_{x \in N_{T}(u)}d_{T}(x)-\sum_{x \in N_{T}(v)}d_{T}(x)\right)+(d_T(x_1)-d_T(x_0))\right]\\
               &>&0.
\end{eqnarray*}
The last inequality follows from $d_T(x_1)>d_T(x_0)$ and $\sum_{x \in N_{T}(u)}d_{T}(x)>\sum_{x \in N_{T}(v)}d_{T}(x).$ Hence, $T \prec_s T'.$
\end{proof}

We are now ready to prove Theorem 1.2.\vspace{2mm}

\noindent {\bf Proof of Theorem 1.2}\ \ Assume that $T$ is the last tree among $\mathscr{T}_\pi$, where $\pi=(d_0,d_1,\ldots,d_{n-1})$ with $d_0 \geqslant d_1 \geqslant \cdots \geqslant d_{n-1}.$  Without loss of generality, we may assume that $V_T=\{v_0,\dots,v_{n-1}\}$ such that $d_T(v_i)\geqslant d_T(v_j)$ for $i<j$, i.e., they are denoted with respect to $d_T(v)$ in non-increasing order. Put $V_i=\{v: {\rm dist}(v,v_0)=i\}$ for $i=0,1,\ldots,p+1$ such that $V_T= \bigcup^{p+1}_{i=0}V_i.$  Denote by $|V_i|=s_i$ for $i=0,1,\ldots,p+1.$ We now may relabel the vertices of $V_T$ by the recursion method. For $V_0$, we relabel $v_0$ by $v_{01}$ and take it as the root of tree $T$. For all vertices of $V_1$ which consists of all neighbors of vertices in $V_0$, may be relabeled as
$$
v_{11}, v_{12}, \ldots, v_{1,s_1}
$$
and satisfy the following conditions:
$$
d_T(v_{11}) \geqslant d_T(v_{12}) \geqslant \cdots \geqslant d_T(v_{1,s_1})
$$
and
\[\label{eq:3.6}
\text{
$
d_T(v_{1i})=d_T(v_{1j})
$\ \
implies\ \
$
\min\{d_{T^*}(x)|x \in N_{T^*}(v_{1i})\cap V_2\} \geqslant \max\{d_{T^*}(x)|x \in N_{T^*}(v_{1j}) \cap V_2\}
$}\]
for $1 \leqslant i < j \leqslant s_1.$ Moreover, $s_1=d_T(v_{01}).$ Generally, we assume that all vertices of $V_i$ are relabeled $\{v_{i1},\dots,v_{i,s_i}\}$ for $i=1,\dots,t.$ Now consider all vertices in $V_{t+1}.$ Since $T$ is a tree, it is easy to see that
$$
s_{t+1}=|V_{t+1}|=d_T(v_{t1})+\cdots+d_T(v_{t,s_t})-s_t.
$$
Hence for $1 \leqslant r \leqslant s_t,$ all neighbors in $V_{t+1}$ of $v_{tr}$ are relabeled as
$$
v_{t+1,d_T(v_{t1})+ \cdots + d_T(v_{t,r-1})-(r-1)+1},\dots, v_{t+1,d_T(v_{t1})+ \cdots + d_T(v_{t,r})-r}
$$
and satisfy the conditions:
$$
d_T(v_{t+1,i})\geqslant d_T(v_{t+1,j}).
$$
and
\[\label{eq:3.7}
\text{$d_T(v_{t+1,i})= d_T(v_{t+1,j})$, implies
$
\min\{d_T(x)|x \in N_T(v_{t+1,i})\cap V_{t+2}\} \geqslant \max\{d_T(x)|x \in N_T(v_{t+1,j}) \cap V_{t+2}\}
$
}\]
for $d_T(v_{t1})+\cdots+d_T(v_{t,r-1})-(r-1)+1 \leqslant i < j \leqslant d_T(v_{t1})+ \cdots + d_T(v_{t,r})-r.$ In this way, we have relabeled all vertices of $V_T=\bigcup_{i=0}^{p+1}V_i.$ Therefore, we are able to define a well ordering of vertices in $V_T$ as follows:
\begin{eqnarray}
v_{ik} \prec v_{jl} \; \text{if}\; 0 \leqslant i < j \leqslant p+1  \; \text{or} \; i=j   \; \text{and} \; 1 \leqslant k < l \leqslant s_i.
\end{eqnarray}
We need to show that this well ordering is a BFS-ordering of tree $T.$ In other words, $T$ is isomorphic to $T^*$.

In order to show this assertion, we first prove that the following  equation holds.
\begin{eqnarray}\label{eq:3.9}
d_{T}(v_{h1}) \geqslant d_{T}(v_{h2}) \geqslant \cdots \geqslant d_{T}(v_{h,s_h}) \geqslant d_{T}(v_{h+1,1})
\end{eqnarray}
for $h=0,1,\ldots,p+1$ by the induction on $h.$

For $h=0,$ clearly, (\ref{eq:3.9}) holds. Assume that (\ref{eq:3.9}) holds for $h=t.$ We consider $h=t+1.$ Suppose to the contrary that $d_{T}(v_{t+1,i}) < d_{T}(v_{t+1,j})$ for $1 \leqslant i < j \leqslant s_{t+1}.$ Then there exist two vertices $v_{tk}$ and $v_{tl}$ with $k < l$ in layer $t$ such that $v_{tk}v_{t+1,j} \in E_T$ and $v_{tl}v_{t+1,j} \in E_T.$  By the induction hypothesis, we have $d_T(v_{tk}) \geqslant d_T(v_{tl}).$ Here we consider the following two possible cases.\vspace{2mm}

{\bf Case 1}.\ $d_T(v_{tk}) > d_T(v_{tl}).$ Let
$$
T'=T-\{v_{tk}v_{t+1,i},v_{tl}v_{t+1,j}\} + \{v_{tk}v_{t+1,j},v_{tl}v_{t+1,i}\}.
$$
It is routine to check that $T' \in \mathscr{T}_\pi.$ By Lemma 3.2, $T \prec_s T',$ a contradiction to $T$ being the last tree among $\mathscr{T}_\pi$.\vspace{2mm}

{\bf Case 2}.\ $d_T(v_{tk}) = d_T(v_{tl}).$ In this case, by (\ref{eq:3.7}) we have that
$\min\{d_{T}(x)|x \in N_{T}(v_{tk})\cap V_{t+1}\} \geqslant \max\{d_{T}(x)|x \in N_{T}(v_{tl}) \cap V_{t+1}\},$ which implies that
\[\label{eq:3.100}
 \sum_{x \in N_T(v_{tk})\cap V_{t+1}}d_T(x) \geqslant \sum_{x \in N_T(v_{tl})\cap V_{t+1}}d_T(x).
\]
Note that there exists a unique vertex, say $v_{t-1, r}$ (resp. $v_{t-1, s}$), in $N_T(v_{tk})$ (resp. $N_T(v_{tl})$) such that $r\le s$. By induction, $d_T(v_{t-1, r})\ge d_T(v_{t-1, s})$, hence together with (\ref{eq:3.100}) we have
$$
 \sum_{x \in N_T(v_{tk})}d_T(x) \geqslant \sum_{x \in N_T(v_{tl})}d_T(x).
$$
Set
\begin{eqnarray}\label{eq:3.110}
  d_T(x_0):= \min\{d_T(x)|x \in N_{T}(v_{tk})\cap V_{t+1}\},\ \ \ \
  d_T(x_1):=\max\{d_T(x)|x \in N_{T}(v_{tl}) \cap V_{t+1}\}.
\end{eqnarray}
Note that $d_{T}(v_{t+1,i}) < d_{T}(v_{t+1,j})$, in view of (\ref{eq:3.110}) we have $d_T(x_0)< d_T(x_1).$
Let
$$
T'=T-\{v_{tk}x_0,v_{tl}x_1\} + \{v_{tk}x_1,v_{tk}x_0\}.
$$
It is routine to check that $T' \in \mathscr{T}_\pi.$ By Lemma 3.3, we have $T \prec_s T',$ a contradiction to $T$ being the last tree among $\mathscr{T}_\pi$.

Similarly, we also show that $d_{T}(v_{h,s_h}) \geqslant d_{T}(v_{h+1,1}).$ Hence (\ref{eq:3.9}) holds.
Therefore, we have
\begin{eqnarray}
d_T(v_{01}) \geqslant d_T(v_{11}) \geqslant \cdots \geqslant d_T(v_{1,s_1}) \geqslant d_T(v_{21}) \geqslant \cdots \geqslant d_T(v_{2,s_2}) \geqslant \cdots \geqslant d_T(v_{p+1,s_{p+1}})
\end{eqnarray}
and
\begin{equation}
\begin{split}
&d_T(v_{01})=d_0,\, d_T(v_{11})=d_1,\, \ldots,\, d_T(v_{1,s_1})=d_{s_1},   \\
&d_T(v_{21})=d_{s_1 +1},\, \ldots,\, d_T(v_{2,s_2})=d_{s_1 + s_2},\,\ldots,   \\
&d_T(v_{p+1,1})=d_{s_1 + \cdots + s_p +1},\, \ldots,\, d_T(v_{p+1,s_{p+1}})=d_{n-1}.
\end{split}
\end{equation}
By (3.8), (3.11) and (3.12), it is easy to see that this well ordering satisfies the conditions (1)-(3) in Definition 2.1. Hence $T$ has a BFS-ordering. Further, by Proposition 1, $T^*$ is isomorphic to $T$. So $T^*$ is the last tree, in an $S$-order, among  $\mathscr{T}_\pi$.
\qed

\begin{lem}
Let $T=(V_T,E_T)$ be a tree with $uv_i \in E_T$ and $wv_i \notin E_T$ for $i=1,2,\dots,k.$ Let $T'=(V_{T'},E_{T'})$ be a new tree from $T$ by deleting edges $uv_i$ and adding edges $wv_i$ for $i=1,2,\dots,k.$ If $d_{T}(w)\geqslant d_{T}(u),$ then $T \prec_s T'.$
\end{lem}
\begin{proof}
Note that $S_i(T)=S_i(T')$ for $i=0,1,2,3,$ by Lemma 2.2(i), we have
\begin{eqnarray*}
S_4(T')-S_4(T)&=&  4(\phi_{T'}(P_3)-\phi_T(P_3))\\
                 &=& 4\left[{d_{T'}(w) \choose 2} + {d_{T'}(u)\choose 2} - {d_{T}(w)\choose 2}-{d_{T}(u)\choose 2}\right]\\
                 &=& 4\left[{d_{T}(w)+k \choose 2} + {d_{T}(u)-k \choose 2} - {d_{T}(w)\choose 2}-{d_{T}(u)\choose 2}\right]\\
                 &=& 4[k^2 + k(d_T(w)-d_T(u))]> 0.
\end{eqnarray*}
Hence, we have $T \prec_s T'.$
\end{proof}

\begin{thm}
Let $\pi$ and $\pi'$ be two different tree degree sequences with the same order. Let $T^*$ and $(T')^*$ be the last tree, in an $S$-order, among $\mathscr{T}_\pi$ and $\mathscr{T}_{\pi'}$, respectively. If $\pi \lhd \pi'$, then $T^* \prec_s (T')^*.$
\end{thm}
\begin{proof}
Note that $\pi \lhd \pi'$, hence by Proposition 2, we may assume, without loss of generality, that $\pi=(d_0,\dots,d_{n-1})$ and $\pi'=(d_0',\dots,d_{n-1}')$ with $d_i=d_i'$ for $i\neq p,\;q,$ and $d_p=d_p'-1,\;d_q=d_q'+1,\; 0\leqslant p < q \leqslant n-1.$ Moreover, let $\pi$ and $\pi'$ be degree sequences of $T^*$ and $(T')^*$, respectively. Since $d_q=d_q'+1 \geqslant 2,$ there exists a vertex $w$ in $N(v_q)$ such that $wv_p \not\in E_{T^*}$ (otherwise, $T^*$ contains a cycle, a contradiction). Let $T_1$ be a tree from $T^*$ by adding the edges $wv_p$ and deleting $wv_q.$ Note that $d_T(v_p)\ge d_T(v_q)$, hence by Lemma 3.4, $T^* \prec_s T_1$. Furthermore, it is easy to see that $T_1\in \mathscr{T}_{\pi'}.$ Hence, $T_1\preceq_s (T')^*$ with equality if and only if $T_1\cong (T')^*$.  Hence $T^* \prec_s T_1 \preceq_s (T')^*,$ as desired.
\end{proof}

From Theorems 1.2 and 3.5, we may deduce the last graph in an $S$-order in some class of graphs. For example, let $\mathscr{T}_{n,s}^{1}$ be the set of all trees of order $n$ with $s$ leaves, $\mathscr{T}_{n,\Delta}^{2}$ be the set of all trees of order $n$ with the largest degree $\Delta$, $\mathscr{T}_{n,\alpha}^{3}$ be the set of all trees of order $n$ with the independence number $\alpha$ and $\mathscr{T}_{n,\beta}^{4}$ be the set of all trees of order $n$ with the matching number $\beta.$

\begin{cor}
A tree $T_1$ is the last tree, in an $S$-order, among $\mathscr{T}_{n,s}^{1}$ if and only if $T_1$ is a star with paths of almost the same length to each of its $s$ leaves (in other words, let $n-1=sq+t,\;0\leqslant t < s$) and $T^*$ is obtained from $t$ paths of order $q+2$ and $s-t$ paths of order $q+1$ by identifying one end of the $s$ paths).

\end{cor}

\begin{proof}
Let $T$ be any tree in $\mathscr{T}_{n,s}^{1}$ with the non-increasing degree sequence $ \pi=(d_0,\dots,d_{n-1})$. Thus $d_{n-s-1} > 1$ and $d_{n-s}=\cdots = d_{n-1}=1$. Let $T^*$ have a BFS ordering tree with the degree sequence $\pi^*=(s,2,\dots,2,1,\dots,1),$ where there are the number $s$ of 1 in $\pi^*.$ By Definition 1, $T^*$ is a star with paths of almost the same length to each of its s leaves. Moreover, it is easy to see that $\pi \lhd \pi^*.$ By Theorem 3.5, the assertion holds.
\end{proof}
\begin{cor}
A tree $T_2$ is the last tree, in a $S$-order, among $\mathscr{T}_{n,\Delta}^{2}$ with $\Delta \geqslant 3$ if and only if $T_2$ is $T^*$ in Theorem $1.2$ with degree sequence $\pi^*$ which is as follow: Denote $l=\left\lceil \log_{(\Delta - 1)} \frac{n(\Delta -2)+2}{\Delta}\right\rceil - 1$ and $n- \frac{\Delta(\Delta -1)^l -2}{\Delta -2}=(\Delta -1)r+q$ for $0 \leqslant q < \Delta -1.$ If $q=0$, put $\pi=(\Delta,\dots,\Delta,1,\dots,1)$ with the number $\frac{\Delta(\Delta -1)^{l-1} -2}{\Delta -2}+r$ of degree $\Delta$. If $1 \leqslant q,$ put $\pi^*=(\Delta,\dots,\Delta,q,1,\dots,1)$ with the number $\frac{\Delta(\Delta -1)^{l-1} -2}{\Delta -2}+r$ of degree $\Delta$.
\end{cor}
\begin{proof}
For any tree $T$ of order $n$ with the largest degree $\Delta,$ let $\pi=(d_0,\dots,d_{n-1})$ be the non-increasing degree sequence of $T$. Assume that $T^*$ has $l+2$ layers. Then there is a vertex in layer 0 (i.e., root), there are $\Delta$ vertices in layer $1$, there are $\Delta(\Delta -1)$ vertices in layer 2,$\dots$, there are $\Delta(\Delta-1)^{l-1}$ vertices in layer $l$, there are at most $\Delta(\Delta -1)^l$ vertices in layer $l+1$. Hence
$$
1+\Delta+\Delta(\Delta-1)+\dots+\Delta(\Delta-1)^{l-1} < n \leqslant 1+\Delta+\Delta(\Delta-1)+\dots+\Delta(\Delta-1)^l.
$$
Thus
$$
\frac{\Delta(\Delta -1)^l -2}{\Delta -2} < n \leqslant \frac{\Delta(\Delta -1)^{l+1} -2}{\Delta -2}.
$$
Hence
$$
l=\left\lceil \log_{(\Delta - 1)} \frac{n(\Delta -2)+2}{\Delta} \right\rceil - 1
$$
and there exist integers $r$ and $0 \leqslant q < \Delta -1$ such that
$$
n- \frac{\Delta(\Delta -1)^l -2}{\Delta -2}=(\Delta -1)r+q.
$$
Therefore degrees of all vertices from layer 0 to layer $l-1$ are $\Delta$ and there are $r$ vertices in layer $l$ with degree $\Delta.$ Denote by $m=\Delta(\Delta -1)^{l-1}-2/\Delta -2 + r-1.$ Then there are $m+1$ vertices with degree $\Delta$ in $T^*$. Hence the degree sequence of $T^* \in \mathscr{T}_{n,\Delta}^{2}$ is $\pi^*=(d_0^*,\dots,d_{n-1}^*)$ with $d_0^*=\cdots=d_m^*=\Delta, \; d_{m+1}^*=\cdots=d_{n-1}^*=1$ for $q=0;$ and is $\pi^*=(d_0^*,\dots,d_{n-1}^*)$ with $d_0^*=\cdots=d_m^*=\Delta, \; d_{m+1}^*=q, \; d_{m+2}^*=\cdots=d_{n-1}^*=1$. It follows from $d_0 \leqslant \Delta$ that $\sum_{i=0}^kd_i \leqslant \sum_{i=0}^kd_i^*$ for $k=0,\dots,m.$ Further by $d_i^*=1 \leqslant d_i$ for $k=m+2,\dots,n-1,$ we have
$$
\sum_{i=0}^kd_i=2(n-1)-\sum_{k+1}^{n-1}d_i \leqslant 2(n-1)- \sum_{k+1}^{n-1}d_i^*=\sum_{i=0}^kd_i^*
$$
for $k=m+1,\dots,n-1$. Thus $\pi \lhd \pi^*.$ Hence by Theorems 1.2 and 3.5, we have $T \preceq_s T^*$ and $"=_s"$ holds if and only if $T=T^*.$
\end{proof}
\begin{cor}
Let $\mathscr{T}_{n,\alpha}^{3}$ be the set of all trees of order $n$ with the independence number $\alpha.$ A tree $T_3$ is the last tree, in an $S$-order, among $\mathscr{T}_{n,\alpha}^{3}$ if and only if $T_3$ is $T^*$ in Theorem $1.2$ with degree sequence $\pi^*=(\alpha,2,\dots,2,1,\dots,1)$ the numbers $n-\alpha-1$ of $2$ and $\alpha$ of 1.
\end{cor}
\begin{proof}
For any tree $T$ of order $n$ with the independence number $\alpha,$ let $I$ be an independent set of $T$ with the independence number $\alpha$ and $\pi=(d_0,\dots,d_{n-1})$ be the non-increasing degree sequence of $T.$ If there exists a pendent vertex $u$ of degree 1 with $u \notin I$, then there exists a vertex $v \in I$ with $uv \in E_T$. Hence $I \bigcup \{u\}\backslash \{v\}$ is an independent set of $T$ with the size $\alpha.$ Therefore, there exists an independent set of $T$ with $\alpha$ which contains all pendent vertices of $T.$ Hence there are at most $\alpha$ pendent vertices. Then $d_{n- \alpha -1} \geqslant 2$ and $\pi \lhd \pi^*.$ Therefore by Theorems 1.2 and 3.5, the assertion holds.
\end{proof}

\begin{cor}
Let $\mathscr{T}_{n,\beta}^{4}$ be the set of all trees of order $n$ with the matching number $\beta.$ A tree $T_4$ is the last tree, in an $S$-order, among $\mathscr{T}_{n,\beta}^{4}$ if and only if $T_4$ is $T^*$ in Theorem $1.2$ with degree sequence $\pi^*=(n-\beta,2,\dots,2,1,\dots,1)$ and the number $n - \beta$ of 1.
\end{cor}
\begin{proof}
For any tree $T$ of order $n$ with the matching number $\beta,$ let $\pi=(d_0,\dots,d_{n-1})$ be the non-increasing degree sequence of $T$. Let $M$ be a matching of $T$ with the matching number $\beta.$ Since $T$ is connected, there are at least $\beta$ vertices in $T$ such that their degrees are at least two. Hence $d_{\beta -1} \geqslant 2.$ Then $\pi \lhd \pi^*.$ Therefore by Theorems 1.2 and 3.5, the assertion holds.
\end{proof}

\end{document}